\documentclass[12pt]{amsart}
\usepackage{amssymb,latexsym}
\usepackage{fullpage}
\usepackage{tikz-cd}
\usepackage[Symbolsmallscale]{upgreek}
\usepackage{array}

\renewcommand{\thefootnote}{}

\begin{document}

\thispagestyle{empty}

\def\theequation {\arabic{section}.\arabic{equation}}
\renewcommand{\thefootnote}{}

\newcommand{\codim}{\mbox{{\rm codim}$\,$}}
\newcommand{\stab}{\mbox{{\rm stab}$\,$}}
\newcommand{\lra}{\mbox{$\longrightarrow$}}

\newcommand{\blr}{\Big \longrightarrow}
\newcommand{\da}{\Big \downarrow}
\newcommand{\ua}{\Big \uparrow}
\newcommand{\hra}{\hookrightarrow}
\newcommand{\rt}{\mbox{\Large{$\rightarrowtail$}}}
\newcommand{\dua}{\begin{array}[t]{c}
\Big\uparrow \\ [-4mm]
\scriptscriptstyle \wedge \end{array}}

\newcommand{\be}{\begin{equation}}
\newcommand{\ee}{\end{equation}}

\newtheorem{guess}{Theorem}[section]
\newcommand{\bth}{\begin{guess}$\!\!\!${\bf }~}
\newcommand{\eeth}{\end{guess}}
\renewcommand{\bar}{\overline}
\newtheorem{propo}[guess]{Proposition}
\newcommand{\bpropo}{\begin{propo}$\!\!\!${\bf }~}
\newcommand{\epropo}{\end{propo}}

\newtheorem{lema}[guess]{Lemma}
\newcommand{\blem}{\begin{lema}$\!\!\!${\bf }~}
\newcommand{\elem}{\end{lema}}

\newtheorem{defe}[guess]{Definition}
\newcommand{\bdefe}{\begin{defe}$\!\!\!${\bf }~}
\newcommand{\edefe}{\end{defe}}

\newtheorem{coro}[guess]{Corollary}
\newcommand{\bcor}{\begin{coro}$\!\!\!${\bf }~}
\newcommand{\ecor}{\end{coro}}

\newtheorem{rema}[guess]{Remark}
\newcommand{\brem}{\begin{rema}$\!\!\!${\bf }~\rm}
\newcommand{\erem}{\end{rema}}

\newtheorem{exam}[guess]{Example}
\newcommand{\beg}{\begin{exam}$\!\!\!${\bf }~\rm}
\newcommand{\eeg}{\end{exam}}

\newcommand{\ctext}[1]{\makebox(0,0){#1}}
\setlength{\unitlength}{0.1mm}
\newcommand{\cl}{{\mathcal L}}
\newcommand{\cn}{{\mathcal N}}
\newcommand{\cp}{{\mathcal P}} 
\newcommand{\ci}{{\mathcal I}}
\newcommand{\bz}{\mathbb{Z}}
\newcommand{\bl}{\mathbb{L}}
\newcommand{\cf}{{\mathcal F}}
\newcommand{\cd}{{\mathcal D}}
\newcommand{\cs}{{\mathcal s}}
\newcommand{\cv}{{\mathcal V}} 
\newcommand{\ce}{{\mathcal E}}
\newcommand{\ck}{{\mathcal K}} 
\newcommand{\ch}{{\mathcal H}}
\newcommand{\cy}{{\mathcal Y}}
\newcommand{\cF}{{\mathcal F}}
\newcommand{\cm}{{\mathcal M}}
\newcommand{\cb}{{\mathcal B}}
\newcommand{\ca}{{\mathcal A}}
\newcommand{\cR}{{\mathcal R}}
 \newcommand{\bq}{\mathbb{Q}}
\newcommand{\bt}{\mathbb{T}} 
\newcommand{\bh}{\mathbb{H}}
\newcommand{\br}{\mathbb{R}} 
\newcommand{\wt}{\widetilde}
\newcommand{\im}{{\rm Im}\,} 
\newcommand{\bc}{\mathbb{C}}
\newcommand{\bp}{\mathbb{P}} 
\newcommand{\bs}{\mathbb{S}} 
\newcommand{\fv}{\mathfrak{v}} 
\newcommand{\fe}{\mathfrak{e}} 
\newcommand{\spin}{{\rm Spin}\,}
\newcommand{\ds}{\displaystyle} 
\newcommand{\tor}{{\rm Tor}\,}

\def\ns{\mathop{\lr}}
\def\nssup{\mathop{\lr\,sup}} 
\def\nsinf{\mathop{\lr\,inf}}
\newcommand{\e}{{\phi}}
\newcommand{\co}{{\cal O}}
\newcommand{\ct}{{\cal T}}

\title{$K$-theory of $\mbox{HYPERK\"AHLER}$ toric manifolds}
\author{V.Uma}     
\address{Department of Mathematics, I.I.T Madras, Chennai-36}

\maketitle

\footnote{2000 A.M.S. Subject Classification:- 14M25, 14F15\\ Key
words and phrases:$\mbox{HyperK\"ahler}$ toric manifolds, hyperplane
arrangements, $K$-theory, $\mbox{hyperK\"ahler}$ toric varieties,
toric quiver varieties}

\noindent {\bf Abstract:} Let $X$ be a toric $\mbox{hyperK\"ahler}$ manifold. The
purpose of this note is to describe the topological $K$-ring $K^*(X)$
of $X$. We give a presentation for the topological $K$-ring in terms
of generators and relations similar to the known
description of the cohomology ring of these manifolds.

\section{Introduction}
Toric $\mbox{hyperK\"ahler}$ manifolds were defined by Bielawski and
Dancer \cite{bd} and have been widely studied recently. The cohomology
ring of a toric $\mbox{hyperK\"ahler}$ manifold has been described by
Konno \cite{konno}. In \cite{str}, Hausel and Sturmfels gave the
algebraic geometric construction of toric $\mbox{hyperK\"ahler}$
varieties and its relation with toric quiver varieties. These
varieties carry the underlying structure of a toric
$\mbox{hyperK\"ahler}$ manifold see \cite[Section 5]{str}.

In this paper we shall use the topological description of toric
$\mbox{hyperK\"ahler}$ manifolds in \cite{konno}. In \cite{konno},
Konno showed the existance of certain canonical complex line bundles
on the $\mbox{hyperK\"ahler}$ manifolds and proved that their first
Chern classes generate the integral cohomology ring. In \cite[Theorem
3.2]{konno}, he further gave the presentation of the cohomology ring
in terms of the combinatorics of certain smooth hyperplane arrangement
naturally associated to the toric $\mbox{hyperK\"ahler}$ manifold.

More recently topological $\mbox{hyperK\"ahler}$ manifolds have been
studied by Kuroki \cite{kur} from the viewpoint of cohomological
rigidity problem, where he has described the equivariant cohomology
ring of these manifolds. In \cite[Theorem 1.1]{str}, Hausel and
Sturmfels gave a presentation of the cohomology ring of toric
$\mbox{hyperK\"ahler}$ varieties with application to the presentation
of the cohomology ring of toric quiver varieties.

Our main aim in this paper is to exploit the known description of the
cohomology ring and its generators to describe the topological
$K$-ring of these manifolds by applying the results in \cite{ah}. Our
methods are similar to those used by Sankaran \cite{s} in the
description of the topological $K$-ring of smooth complete toric
varieties and torus manifolds. We recall here that the
$\mbox{hyperK\"ahler}$ manifold although non-compact in general (e.g
$T^*(\bp_{\bc}^n$)) is homotopy equivalent to a finite CW complex
namely its core and hence has the structure of a CW complex of finite
type. The key tool is the application of the Atiyah-Hirzebruch
spectral sequence which degenerates in this setting.  Using this in
Theorem \ref{main} we show that the isomorphism classes of the
canonical line bundles on the $\mbox{hyperK\"ahler}$ manifold defined
by Konno generate the topological $K$-ring. The presentation of the
$K$-ring follows from that of the cohomology ring as in the case of
toric manifolds.

One of the motivations of this paper is to understand the topological
$K$-ring of toric $\mbox{hyperK\"ahler}$ varieties and toric quiver
varieties developed by Hausel and Sturmfels in \cite{str}. Since toric
$\mbox{hyperK\"ahler}$ varieties are biholomorphic to toric
$\mbox{hyperK\"ahler}$ manifolds, the presentation of their
topological $K$-rings will follow from our main theorem. However
because of their additional algebraic geometric and combinatorial
structure, we expect that the presentation of their $K$-ring gets a
canonical interpretation in terms of matroid ideals and circuit
ideals, similar to that of their cohomology ring described in
\cite[Theorem 1.1]{str}. This is work in progress.

Let $T'$ be the compact torus of one fourth the dimension of the
$\mbox{hyperK\"ahler}$ toric manifold which acts canonically on it.
We believe that we can give a description of the $T'$ equivariant
$K$-ring of these manifolds using its $GKM$ structure similar to the
description of the equivariant cohomology ring by Kuroki in
\cite[Section 5]{kur}. This shall be taken up in future work. 

There is also the simplicial analogue of $\mbox{hyperK\"ahler}$ toric
orbifolds (see \cite{bd},\cite{konno}, \cite{str} and
\cite{kur}). They are associated to simple and not necessarily smooth
hyperplane arrangement.  In this paper we shall unless otherwise
specified work with smooth $\mbox{hyperK\"ahler}$ manifolds.

Let $X$ denote a $4n$ dimensional $\mbox{hyperK\"ahler}$ toric
manifold equipped with the action of the $n$-dimensional compact torus
$T':=(S^1)^n$ preserving the $\mbox{hyperK\"ahler}$ structure. Let
$H_i$, $1\leq i\leq m$ denote the associated smooth hyperplane
arrangement in the dual of the Lie algebra
$\mathfrak{t}'^*\simeq \mathbb{R}^n$ of $T'$. By \cite{konno} we know
that there exists $m$ complex line bundles $L_i$, $1\leq i\leq m$ on
$X$ such that $\tau_i:=c_1(L_i)\in H^2(X;\mathbb{Z})$ generate
$H^*(X;\mathbb{Z})$. Let $\mathfrak{J}'$ denote the ideal in
$\bz[x_1,\ldots,x_m]$ defined by the following two sets of relations:
\be\label{preshyperk1} \prod_{i\in I}x_i ~~\mbox{whenever}~~ \bigcap_{i\in
  I}H_i=\emptyset ~\mbox{for}~I\subseteq [1,m]\ee \be\label{preshyperk2}
\prod_{j, \langle u,v_j \rangle>0} (1-x_j)^{\langle u,v_j
  \rangle}-\prod_{j, \langle u,v_j \rangle<0} (1-x_j)^{-\langle
  u,v_j\rangle} ~\mbox{for}~ u\in \mathfrak{t}_{\bz}'^*.\ee The following
is our main theorem which describes the topological complex $K$-ring
of $X$.

{\bf Main Theorem:} The map from $\bz[x_1,\ldots,x_m]$ to $K^*(X)$
which sends $x_j\mapsto 1-[L_j]$ defines a ring homomorphism
$\psi:\mathcal{R}:=\bz[x_1,\ldots,x_m]/\mathfrak{J}'\lra
K^*(X)$. Moreover, $\psi$ defines an isomorphism of $\bz$-algebras.

\section{Definition of toric $\mbox{hyperK\"ahler}$ manifolds}

We begin by briefly recalling the definition of $\mbox{hyperK\"ahler}$
manifold, and required terminologies and notations on (see
\cite[Section 3]{bd} \cite{konno}, \cite{kur} and \cite[Section
5]{str}).
 
Recall that multiplication by $i$ (resp. $j$ and $k$) defines three
complex structures $I$ (resp. $J$ and $K$) on the $m$ dimensional
quaternionic vector space $\bh^m$ which satisfy the quaternionic
relations. Consider the Euclidean metric $g$ on
$\bh^m\simeq \br^{4m}\simeq \br^m\oplus i\br^m\oplus j\br^m\oplus
k\br^m$. We define the $\mbox{K\"ahler}$ forms on $\bh^m$
\be\label{kahler} \omega_I(X,Y)=g(IX,Y)\ee (and similarly
$\omega_J, \omega_K$) where $X$ and
$Y$ are tangent vectors at a point in $\bh^m$. Then
$g$ is a {\it
  $\mbox{hyperK\"ahler}$ metric} that is a
$\mbox{K\"ahler}$ metric with respect to all the three complex
structures.

The symplectic group $Sp(m)\subseteq SO(4m)$ consists of matrices
which commute with $I$, $J$ and $K$. Then $Sp(m)$ preserves the
$\mbox{hyperK\"ahler}$ metric or equivalently preserves the
$\mbox{K\"ahler}$ forms $\omega_I,\omega_J$ and $\omega_K$. We
consider the action of $Sp(m)$ on $\bh^m$ from the right.

We fix the identification ${\bh}^m\rightarrow\bc^m\times \bc^m$ such
that
$\xi:=(\xi_1,\ldots,\xi_m)\mapsto ((z,w):=(z_1,\ldots,z_m), (w_1,\ldots, w_m))$
where $\xi_r=z_r+w_r J$ for $z_j,w_j\in \br+\br\cdot I\simeq \bc$
and $1\leq r\leq m$.

The diagonal subgroup $T=(S^1)^m\subseteq Sp(m)$  acts on $\bh^m$
as
\[e^{i\theta}\cdot (z_1+w_1J,\ldots, z_m+w_mJ)=((z_1+w_1J)\cdot e^{i
    \cdot\theta_1},\ldots, (z_m+w_mJ)\cdot e^{i \cdot\theta_m})\]
where
$e^{i\theta}:=(e^{ i\cdot\theta_1},\ldots, e^{ i
  \cdot\theta_m})\in T$. Using the quaternionic relation
$i\cdot j=-j\cdot i=k$ the action of $T$ on $\bh^m$ can be rewritten
as follows:
$\xi\cdot e^{i\theta}:=(z\cdot e^{i\theta}, w\cdot e^{-i\theta})$.

Let $\mathfrak{t}$ denote the Lie algebra of $T$.  Since the action
preserves the $\mbox{hyperK\"ahler}$ structure, it gives the
$\mbox{hyperK\"ahler}$ moment map
\[\mu=(\mu_{I}, \mu_{J}, \mu_{K})=(\mu_{\mathbb{R}},\mu_{\mathbb{C}}):\bh^m\rightarrow
  \mathfrak{t}^*\otimes \mathbb{R}^3.\] where
$\mu_{I}, \mu_{J}, \mu_{K}$ are the
$\mbox{K\"ahler}$ moment maps with respect to
$\omega_I, \omega_J, \omega_K$ respectively.

Let $\{e_1,\ldots, e_m\}$ denote the standard basis of
$\mathfrak{t}\simeq \mathbb{R}^m$ and
$\{e^*_1,\ldots, e^*_m\}\subseteq \mathfrak{t}^*$ be the dual basis. Let
$\displaystyle{\mathfrak{t}_{\bz}:=\sum_{r=1}^m\bz e_r,~~~
\mathfrak{t}_{\bz}^*:=\sum_{r=1}^m\bz e^*_r}$ We have the
$\mbox{K\"ahler}$ moment maps $\mu_{I},\mu_J,\mu_K$ with respect to
$\omega_I, \omega_J,\omega_K$ respectively. We can write
\[\mu_{\br}:=\mu_{I}(z,w)=\frac{1}{2}\sum_{r=1}^m
  (|z_r|^2-|w_r|^2)e^*_r\in \mathfrak{t}^*\]
\[\mu_{\bc}:=
  (\mu_{J}+i\mu_{K})(z,w)=\sum_{r=1}^mz_rw_re^*_r\in
  \mathfrak{t}^*\otimes \bc=\mathfrak{t}_{\bc}^*.\] Here $\mu_{\bc}$
is the moment map for the $I$-holomorphic action of
$T_{\bc}\simeq (\bc^*)^m$ on $\bh^m$ with respect to the holomorphic
symplectic form $\omega_{\bc}=\omega_{J}+i\omega_K$.

Let $H\simeq (S^1)^{m-n}$ denote the subtorus of $T$ which is the kernel of
a surjective map of tori $T\lra T'$ where $T'\simeq T/H$ is a torus of
dimension $n$. We therefore get the following short exact sequence of
compact tori \be\label{es3} 0\longrightarrow
H\stackrel{{\iota}}{\longrightarrow}
T\stackrel{{\rho}}{\longrightarrow} T'\longrightarrow 0\ee where
$\iota:H\hra T$ is the inclusion and $\rho:T\lra T'$ is the projection.

We also have the following short exact sequences of the corresponding Lie
algebras and their duals \be\label{es1} 0\longrightarrow
\mathfrak{h}\stackrel{\iota_{*}}{\longrightarrow}
\mathfrak{t}\stackrel{\rho_{*}}{\longrightarrow}
\mathfrak{t}'\longrightarrow 0\ee \be\label{es2} 0\longrightarrow
\mathfrak{t}'^*\stackrel{\rho^*}\longrightarrow \mathfrak{t}^*
\stackrel{\iota^*}\longrightarrow \mathfrak{h}^*\longrightarrow 0.\ee
Let $\iota_{\bc}^*: \mathfrak{t}_{\mathbb{C}}^*\lra
\mathfrak{h}_{\mathbb{C}}^*$ denote the map induced on the
complexification of the dual of the Lie algebras

Since $H$ acts on $\bh^m$ preserving its $\mbox{hyperK\"ahler}$
structure, we have the $\mbox{hyperK\"ahler}$ moment map
$\mu_{H}:=(\iota^*\oplus\iota_{\bc}^*)\circ (\mu_{\br}\oplus
\mu_{\bc}):\bh^m\rightarrow \mathfrak{h}^*\oplus
\mathfrak{h}^*_{\mathbb{C}}$.

For $\alpha\neq 0$ in $\mathfrak{h}^*$,
$(\alpha,0)\in \mathfrak{h}^*\oplus \mathfrak{h}^*_{\bc}$ is a regular
value of $\mu_{H}$. The necessary and sufficient condition for the
action of $H$ on $\mu_{H}^{-1}(\alpha,0)$ to be free or equivalently
for $\mu_{H}^{-1}(\alpha,0)/H$ to be a smooth manifold is that we have
the following split exact sequence of lattices \be\label{esint}
0\longrightarrow \mathfrak{h}_{\bz}\stackrel{\iota_*}\longrightarrow
\mathfrak{t}_{\bz} \stackrel{\rho_*}\longrightarrow
\mathfrak{t}'_{\bz}\longrightarrow 0\ee (see \cite[Proposition
2.2]{konno},\cite{bd} and \cite[Proposition 2.1]{kur}). We call the
manifold $X:=X(\alpha,0)=\mu_{H}^{-1}(\alpha,0)/H$ of real dimension
$4n$ a {\it toric $\mbox{hyperK\"ahler}$ manifold}. Since $\iota^*$ is
surjective we fix an element $\text{v}\in \mathfrak{t}^*$ such that
$\iota^*(\text{v})=\alpha$. Moreover, by \cite[Theorem 2.1]{konno} the
diffeomorphism type of $X$ is independent of the regular value chosen.

\beg\label{cotproj} Let $T=(S^1)^{n+1}$, $T'=(S^1)^{n}$ and
$\rho:T\longrightarrow T'$ be a map of tori defined by
$(t_0,\ldots,t_{n})\mapsto
(t_0t_{n}^{-1},\ldots,t_{n-1}t_{n}^{-1})$. The subgroup
$H=\mbox{ker}(\rho)\simeq S^1$ is the diagonal subgroup
$(t,t,\ldots, t)\in T$.  The complex projective space $\bc\bp^n$ is
the $\mbox{K\"ahler}$ quotient of $\bc^{n+1}$ by the $H$ action given
by scalar multiplication and the $\mbox{K\"ahler}$ form
$\displaystyle{\frac{i}{2}\sum_{r=1}^{n+1}dz_r\wedge d\bar{z}_r}$.  On
$\bh^{n+1}=\bc^{n+1}+\bar{\bc}^{n+1}=\bc^{n+1}\bigoplus
(\bc^{n+1})^{*}$ we consider the action of $H$ given by
$(z,w)\cdot t:=(z\cdot t,w\cdot t^{-1})$. The $\mbox{K\"ahler}$ form
for the complex structure on $\bh^{n+1}$ induced from the complex
structure on $\bc^{n+1}$ is given by
$\displaystyle{\omega_{I}:=
\frac{i}{2}\sum_{r=1}^{n+1}}dz_r\wedge
d\bar{z}_r+dw_r\wedge d\bar{w}_r$. Also
$\displaystyle{\omega_{\bc}=\omega_{J}+i\omega_K=\sum_{r=1}^{n+1}dz_r\wedge
  dw_r}$ is the corresponding holomorphic symplectic form. The moment
map $\mu_{H}$ is given by
$\mu_{H}(z,w)=(\frac{1}{2}(|z|^2-|w|^2),  z\cdot w)$. Taking the
regular value $(\frac{1}{2},0)\in \br\oplus \bc$ the
$\mbox{hyperK\"ahler}$ quotient $\mu_{H}^{-1}(\frac{1}{2},0)/H$ is the space
\[\{(z,w)\in \bc^{n+1}\oplus (\bc^{n+1})^* ~:~z\neq 0, z\cdot w=0\}.\]
This can be identified with the cotangent space of $\bc\bp^{n}$ namely
$T^*_{[z]}\bc\bp^{n}:=[z]\otimes (\bc^{n+1}/[z])^*$.\eeg

\subsection{Hyperplane arrangement and 
  Line bundles on toric $\mbox{hyperK\"ahler}$ manifolds}

For $i=1,\ldots, m$, let
$H_i:=\{ p\in \mathfrak{t}'^*|\langle \rho^*p+\text{v},
e_i\rangle=0\}$. Then each $H_i$ as an affine hyperplane in
$\mathfrak{t}'^*$ with a choice of nonzero oriented vector
$v_i:=\rho(e_i)$ for $i=1,\ldots, m$ (see \cite[Remark
2.5]{kur}). Equivalently
$H_i:=\{p\in \mathfrak{t}'^*\mid \langle p,v_i\rangle=-\langle
\text{v},e_i\rangle\}.$

Let $\ch:=\{H_1,\ldots,H_m\}$ denote the hyperplane arrangement
associated with the toric $\mbox{hyperK\"ahler}$ manifold
$X(\alpha,0)$. Then $\ch$ is a {\it smooth} hyperplane arrangement
that is (i) whenever $H_{i_1}\cap\cdots\cap H_{i_l}\neq \emptyset$ the
codimension of $H_{i_1}\cap\cdots\cap H_{i_l}$ is $l$ and (ii) the
vectors $v_{i_1},\ldots, v_{i_l}$ can be extended to a $\bz$ basis of
$\mathfrak{t}'_{\bz}$ (see \cite[Section 3]{bd}).

We now recall the definition of characteristic submanifolds of the
toric $\mbox{hyperK\"ahler}$ manifold from \cite[Section 3.3]{kur}.

Let $X_i$ denote the invariant connected submanifold of $X(\alpha,0)$
which is fixed by the circle subgroup $T_i$ of $T'$ obtained by the
exponent of $v_i\in \mathfrak{t}'$. The circle subgroup is determined
by the primitive vector $v_i$ upto sign. Then $T_i$ acts on the normal
bundle $N_i$ of $X_i$ by right scalar multiplication on the fibers
which are isomorphic to $\bh=\bc\oplus\bar{\bc}$. Fixing a sign of
$v_i$ or equivalently an orientation of $N_i$ is called an
omniorientation of $X_i$ for $1\leq i\leq m$.

We can alternately construct the characteristic submanifold $X_i$ as
$\mbox{hyperK\"ahler}$ quotient of the restricted action of $H$ on the
$m-1$ dimensional subspace $Y_i$ of $Y=\bh^m$ defined by the vanishing
of the $i$th coordinate vector (see \cite[Proposition 3.6]{kur}).
Thus the real dimension of $X_i$ is $4n-4$ for $1\leq i\leq m$. The
total space of the normal bundle $N_i$ of $X_i$ is
$E(N_i)=Y_i\times_{H} \bh_i$ by \cite[Proposition 3.6]{kur}. Here
$\bh_i\simeq \bh\simeq \br^4$ denotes the $1$-dimensional $\bh$-vector
space which is the representation of $H$ where $H$ acts by
$\iota_i:=p_i\circ \iota$, corresponding to the $i$th projection $p_i$
of $T$ onto its coordinate circle subgroup $S_i\simeq S^1$. Thus the
normal bundle is the restriction to $X_i$ of the associated bundle
$\mu_{H}^{-1}(\alpha,0)\times_{H} \bh_i$ on $X(\alpha,0)$. By
construction there is a canonical induced action of the torus $T'=T/H$
on $\mu_{H}^{-1}(\alpha,0)\times_{H} \bh_i$ which gives it the
structure of a $T'$-equivariant bundle on $X(\alpha,0)$ (see
\cite[Section 5.2]{kur} for details).

Now, $\bh=\bc\oplus \bar{\bc}$ so that the bundle
$\mu_{H}^{-1}(\alpha,0)\times_{H} \bh_i$ splits into the following
line bundles
$\mu_{H}^{-1}(\alpha,0)\times_{H} (\bc_i\oplus \bar{\bc_i})$ where
$\bc_i$ denotes the $1$-dimensional complex representation of $H$
where $H$ acts on $\bc$ via the character
$\iota_i:=p_i\circ \iota$. Let
\[L_i:=\mu_{H}^{-1}(\alpha,0)\times_{H} \bc_i\] denote the associated
complex line bundle on $X(\alpha,0)$.  We recall from \cite[Section
4]{konno} and \cite[Section 5.2]{kur} that these line bundles are
holomorphic with respect to the complex structure $I$ on $X(\alpha,0)$
and $\tau_i:=c_1(L_i)$ for $1\leq i\leq m$ generate the integral
cohomology ring of $X(\alpha,0)$.

\subsection{Cohomology ring of toric hyperKahler manifolds}

We recall below the presentation of the ordinary integral cohomology
ring of the toric $\mbox{hyperK\"ahler}$ manifold $X$
due to Konno.

\bth\label{cohom} (\cite[Theorem 3.2]{konno}).  Let
$X:=(X(\alpha,0),T')$ be a toric $\mbox{hyperK\"ahler}$ manifold and
$\ch:=\{H_1,\ldots, H_m\}$ be the associated smooth hyperplane
arrangement. Let $\mathfrak{J}$ denote the ideal in
$\bz[x_1,\ldots, x_m]$ defined by the following two sets of relations:
\be\label{preshypereco1} \prod_{i\in I}x_i ~~\mbox{whenever}~~
\bigcap_{i\in I}H_i=\emptyset \ee for $I\subseteq [1,m]$, and
\be\label{preshypereco2}\sum_{j=1}^m\langle u,v_j\rangle x_j\ee for
$u\in \mathfrak{t}_{\bz}'^*$ . Then the canonical map from
$R:=\bz[x_1,\ldots, x_m]/{\mathfrak{J}}$ to $H^*(X(\alpha,0);\bz)$
which sends $x_j\mapsto c_1(L_j)$ is an isomorphism of $\bz$-algebras.
\eeth

\section{Main theorem}

We now state and prove the main theorem of this paper which gives a
presentation for the topological $K$-ring of the toric
$\mbox{hyperK\"ahler}$ manifold $X=X(\alpha,0)$.

\bth\label{main} Let $X:=(X(\alpha,0),T')$ be a toric
$\mbox{hyperK\"ahler}$ manifold and $\ch:=\{H_1,\ldots, H_m\}$ be the
associated smooth hyperplane arrangement. Let
$\mathfrak{J}'$ denote the ideal in $\bz[x_1,\ldots,x_m]$ defined by
the following two sets of relations: \be\label{preshyperk1}
\prod_{i\in I}x_i ~~{whenever}~~ \bigcap_{i\in I}H_i=\emptyset
~{for}~I\subseteq [1,m]\ee \be\label{preshyperk2} \prod_{j, \langle
  u,v_j \rangle>0} (1-x_j)^{\langle u,v_j \rangle}-\prod_{j, \langle
  u,v_j \rangle<0} (1-x_j)^{-\langle u,v_j\rangle} ~\mbox{for}~ u\in
\mathfrak{t}_{\bz}^*.\ee Then the canonical map $\psi$ from
$\mathcal{R}:=\bz[x_1,\ldots,x_m]/\mathfrak{J}'$ to the topological
$K$-ring $K^*(X)$ of $X$ which sends $x_j\mapsto 1-[L_j]$ is an
isomorphism of $\bz$-algebras.  \eeth
 
{\bf Proof:} The proof of this theorem follows along the lines of the
proof \cite[Theorem 2.2]{s}. The relation \ref{preshypereco1} in
$H^*(X(\alpha,0);\bz)$ implies that if
\(\displaystyle H_{i_1}\cap\cdots\cap H_{i_k}=\emptyset\) then
\[c_{k}(L_{i_1}\oplus L_{i_2}\oplus\cdots L_{i_k})=c_1(L_{i_1})\cdot
c_1(L_{i_2})\cdots c_1(L_{i_k})=0.\] Now, applying the
$\gamma$-operation in $K$-theory (see \cite[Proposition 7.4]{karoubi})
we get
\[\gamma^k([L_{i_1}\oplus \cdots \oplus L_{i_k}]-k)=(-1)^k
c_k(L_{i_1}\oplus \cdots \oplus L_{i_k})=0.\] Further, we have
\[\gamma^k([L_{i_1}\oplus \cdots \oplus L_{i_k}]-k)=\left([L_{i_1}]-1
\right)\cdot \left([L_{i_2}]-1 \right)\cdots \left([L_{i_k}]-1
\right).\] This proves the relation (\ref{preshyperk1}) in $K^*(X)$.
For $u\in \mathfrak{t}_{\bz}'^*$ let
\[L_u:=\prod_{i=1}^m L_i^{\langle u, v_i
  \rangle}.\] Then
\(\displaystyle c_1(L_u)=\sum_{i=1}^m\langle u,v_i\rangle c_1(L_i)=0\)
by (\ref{preshypereco2}). This implies that $L_u$ is isomorphic to a
trivial line bundle on $X$. Thus in $K^*(X)$ we have 
\[\prod_{i=1}^N [L_i]^{\langle u,v_i \rangle} =1.\] proving relation
(\ref{preshyperk2}).  Hence $\psi$ defines a ring homomorphism
$\mathcal{R}\longrightarrow K^*(X)$. Further, by Theorem \ref{cohom},
$c_1(L_i), 1\leq i\leq N$ generate the cohomology of $X$.  In
particular, $H^*(X;\bz)$ is generated by $H^2(X;\bz)$ and vanishes in
odd dimensions. Also $X$ has the homotopy type of a finite CW complex
since by \cite[Section 6]{konno}, the $\mbox{Core}(X)$, which is a
finite union of compact toric submanifolds of $X$ is a
$T'$-equivariant deformation retract of $X$. We now proceed by the
arguments similar to \cite[Section 3]{s} and \cite[Lemma 4.1, Lemma
4.2]{su}. Let $f_i:X\longrightarrow \bp^{\infty}$ be a classifying map
for the bundle $L_i$ for $1\leq i\leq m$. We consider the map
$f:X\longrightarrow (\bp^{\infty})^m$ defined as
$f(x)=(f_1(x),\ldots, f_m(x))$. Since $c_1(L_i), 1\leq i\leq m$
generate $H^*(X;\bz)$ it follows that
$f^*:H^*((\bp^{\infty})^m;\bz)\longrightarrow H^*( X;\bz)$ is
surjective. By the naturality of the Atiyah-Hirzebruch spectral
sequence it follows that $K^*(X)$ is generated by $[L_i]$,
$1\leq i\leq m$. Thus we conclude that $\psi$ is
surjective. Furthermore, it can be shown that $R$ and hence
$H^*(X;\bz)$ is a free abelian group of finite rank (see \cite[Theorem
3.2]{konno} and \cite[{\bf 3.8}]{d}) which is equal to the Euler
characteristic $\chi(X)$.  Thus the collapsing of the Atiyah
Hirzebruch spectral sequence implies that $K^*(X)$ is also a free
abelian group of rank equal to $\chi(X)$ (see \cite[p. 19]{ah}).

We shall now show that $\psi$ is injective by proving that
$\mathcal{R}$ is free abelian of rank $\chi(X)$.

As in \cite[Section 4]{s} we construct a filtered ring
$S=\bz[x_1,\ldots, x_m]$ and let $\mathfrak{J}'$ the ideal defined by
the relations (\ref{preshyperk1}) and (\ref{preshyperk2}). The ring
$S$ is graded with $\mbox{deg}(x_i)=1$. The abelian group of all
homogeneous polynomials of degree $j$ is denoted by $S_{(j)}$. We then
have a multiplicative filtration
$S=S_0\supseteq S_1\supseteq \cdots \supseteq S_r\supseteq\cdots$
where \(\displaystyle S_r=\bigoplus_{j\geq r} S_{(j)}\). Since
$\mathfrak{J}'$ is an ideal generated by elements with constant term
zero, the above filtration is $\mathfrak{J}'$-stable. This induces a
decreasing multiplicative filtration
$\cR=\cR_0\supseteq \cR_1\supseteq \cR_2\supseteq \cdots$ on
$\mathcal{R}:=S/\mathfrak{J}'$. Let $\mbox{gr}(\mathcal{R})$ denote
the associated graded ring with respect to this filtration. For
$f\in S_r$ we denote by $\mbox{in}(f)$ the initial form of $f$ which
is the homogeneous polynomial of degree $r\geq 1$ such that
$f-\mbox{in}(f)\in S_{r+1}$.

If
\(\displaystyle z_u:= \prod_{j, \langle u,v_j \rangle>0}
(1-x_j)^{\langle u,v_j \rangle}-\prod_{j, \langle u,v_j \rangle<0}
(1-x_j)^{-\langle u,v_j\rangle}\) for each
$u\in \mathfrak{t}_{\bz}'^*$. Then
\(\displaystyle h_u:=\mbox{in}(z_u)=\sum_{i=1}^N \langle u, v_i\rangle
x_i\). Thus the ideal $\mathfrak{J}$ is
generated by the relations (\ref{preshypereco1}) and $h_u$ for
$u\in \mathfrak{t}_{\bz}'^*$. Since $\mathfrak{J}$ is a graded ideal
in $S$, $R=S/\mathfrak{J}$ is a graded abelian group. Moreover, by the
arguments in \cite[pages 462-463]{s} we have a surjective homomorphism
of graded abelian groups
$\eta: \mbox{gr}(R)=R\longrightarrow \mbox{gr}(\mathcal{R})$. Since
$R$ is free abelian of rank $\chi(X)$ it follows that
$\mbox{gr}(\mathcal{R})$ is free abelian of rank at most
$\chi(X)$. Since $\mathcal{R}$ and $\mbox{gr}(\mathcal{R})$ are free
abelian of the same rank it follows that $\mathcal{R}$ is free abelian
of rank at most $\chi(X)$.  Combining with the surjectivity of $\psi$
this implies that $\psi$ must be an isomorphism.

We now illustrate Theorem \ref{main} by describing the $K$-ring of the
cotangent bundle of the complex projective space whose construction as
a toric $\mbox{hyperK\"ahler}$ manifold was described in Example
\ref{cotproj}. (See \cite[Example 2.4 and Example 3.4]{kur})

\beg\label{cotangent bundle of the complex projective space} From
Example \ref{cotproj} we have the inclusion $H\subset T=(S^1)^{n+1}$
of the diagonal and $T'=T/H\simeq (S^1)^n$. This induces the following
exact sequence of dual of the corresponding Lie algebras
\[\mathfrak{t}'^*\stackrel{\rho^*}{\longrightarrow}\mathfrak{t}^*\stackrel{\iota^*}{\longrightarrow}
  \mathfrak{h}^*\] where
$\iota^*(\alpha_1,\ldots,
\alpha_{n+1})=\alpha_1+\cdots+\alpha_{n+1}\in \br=\mathfrak{h}^*$ where
$(\alpha_1,\ldots,\alpha_{n+1})\in \mathfrak{t}^*\simeq
\br^{n+1}$. Further,
\[\begin{split}\rho^*(u)&=(\langle u,e_1\rangle,\ldots, \langle u,e_n\rangle\langle
u, -(e_1+\ldots+e_n)\rangle)\\ &=(a_1,\ldots, a_n,-a_1-a_2-\cdots -a_n)\end{split}\]
for $u=(a_1,\ldots,a_n)\in \mathfrak{t}'^*\simeq \br^{n}$.  For
$\alpha=n+1$ we choose the lift $(1,1,\ldots, 1)=\text{v}\in
\mathfrak{t}^*$. Thus the associated hyperplane arrangement is given
by
\[\begin{array}{llllll} H_1&=&\{(a_1,\ldots, a_{n})\in
    \mathfrak{t}'^*\mid a_1=-1\}\\ &\vdots& \\H_n&=& \{(a_1,\ldots,
    a_{n})\in \mathfrak{t}'^*\mid a_n=-1\} \\H_{n+1}&=&\{(a_1,\ldots, a_{n})\in
    \mathfrak{t}'^*\mid a_1+a_2+\cdots+a_n=1\} \\ \end{array}.\]

Clearly $I=[1,n+1]$ is the only subset of $[1,n+1]$ such that
$\displaystyle{\bigcap_{i\in I} H_i=\emptyset}$. Hence by Theorem
\ref{main} it follows that the topological $K$-ring of $T^*(\bc\bp^n)$
is isomorphic to $\mathcal{R}:=\bz[x_1,\ldots,x_{n+1}]/\mathfrak{J}'$ where
$\mathfrak{J}'$ is the ideal in $\bz[x_1,\ldots,x_{n+1}]$ generated by
the monomial $x_1\cdot x_2\cdots x_{n+1}$ and the following $n$
relations
\[(1-x_1)-(1-x_{n+1}), (1-x_2)-(1-x_{n+1}),\ldots,
  (1-x_n)-(1-x_{n+1})\] corresponding to the basis
$\{(1,0,\ldots,0), (0,1,\ldots, 0),\ldots, (0,\ldots,1)\}$ of
$\mathfrak{t}'^*$.  Let $L_{n+1}$ be the canonical line bundle on
$T^*(\bc\bp^n)$ corresponding to the hyperplane $H_{n+1}$. After
suitable change of variables in the ring $\mathcal{R}$ it can be seen
that the map $x\mapsto [L_{n+1}]$ defines an isomorphism of
$\bz$-algebras from $\bz[x]/(1-x)^{n+1}$ to $K^*(X)$. \eeg

\noindent {\bf Acknowledgements:} I wish to thank Prof. P. Sankaran
for valuable discussions during the preparation of this manuscript. I
am grateful to him for patiently reading the preliminary version and
his suggestions for improving the presentation.

\end{document}